\documentclass{amsart}

\newtheorem{theorem}{Theorem}[section]
\newtheorem{lemma}[theorem]{Lemma}

\theoremstyle{definition}
\newtheorem{definition}[theorem]{Definition}
\newtheorem{example}[theorem]{Example}

\theoremstyle{remark}

\numberwithin{equation}{section}

\begin{document}

\title{Hyperholomorphic functions on commutative algebras}

\author{Anatoliy A. Pogorui}
\address{Department of Mathematics, Zhytomyr State University, Zhytomyr, Ukraine 10008}
\email{pogor@zu.edu.ua}



\date{November 2, 2006.}



\begin{abstract}
In this paper we study properties of hyperholomorphic functions on
commutative finite algebras. It is investigated the Cauchy-Riemann
type conditions for hyperholomorphic functions. We prove that a
hyperholomorphic function on a commutative finite algebra can be
expanded in a Taylor series.  We also present a technique for
computing zeros of polynomials in commutative algebras.
\end{abstract}

\maketitle

\section*{Introduction}

The development of hyperholomorphic function analysis has renewed
interest in mathematics and physics because of fruitful
applications. One of the most popular hypercompex analysis is
quaternionic analysis, however, non-commutativity of quaternion
algebra causes many intractable problems, for instance, the
problem of expansion of a hyperholomorphic quaternionic function
in a Taylor series. In this regard hyperholomorphic analysis on
commutative unitary algebras is a natural extension of complex
analysis, despite the fact that in these algebras we have problem
of zero divisors. There are many commutative generalization of
complex numbers, say, hyperbolic numbers, bicomplex algebra etc.
(see \cite{1}). In \cite{2} it is proved that hyperholomorphic
functions on bicomplex algebra can be expanded in a Taylor series.
In this article we generalize this result to any finite
commutative algebra.

\section{Differentiation in finite commutative algebra}
Let  $\mathbf{A}$  be a finite commutative unitary algebra over
$K=\mathbb{R}$ (or $\mathbb{C}$), a set of vectors $\vec{e}_{0},
\vec{e}_{1}, \ldots, \vec{e}_{n}$ be a basis of $\mathbf{A}$, and
$\vec{e}_{0}$ be the unit of the algebra.
 Consider a function
$\vec{f}:\mathbf{A}\rightarrow\mathbf{A}$ of the following form
$$
 \vec{f}(\vec{x})=\sum_{k=0}^{n}\vec{e}_{k}u_{k}(\vec{x}),
$$
where $u_{k}(\vec{x})=u_{k}(x_{1},x_{2},\ldots,x_{n})$  are real
(or complex) functions of $n$ arguments.

\begin{definition}
$\vec{f}(\vec{x})$ is called $\mathbf{A}$-differentiable at a
point $\vec{x}_{0}\in \mathbf{A}$ if there exists the function
$\vec{f'}:\mathbf{A}\rightarrow\mathbf{A}$  such that for any
$\vec{h}\in \mathbf{A}$
\begin{equation}\label{1}
\vec{h}\vec{f'}(\vec{x}_{0})=\lim_{\varepsilon\rightarrow 0}
\frac{\vec{f}(\vec{x}_{0}+\varepsilon\vec{h})-\vec{f}(\vec{x}_{0})}{\varepsilon},
\end{equation}
where $\vec{f'}$ doesn't depend on $\vec{h}.$
\end{definition}
A function $\vec{f}$  is said to be $\mathbf{A}$-holomorphic if
$\vec{f}$ is $\mathbf{A}$-differentiable at every point of
$\mathbf{A}$.
\begin{theorem} \label{A}
A function
$\vec{f}(\vec{x})=\sum_{k=0}^{n}\vec{e}_{k}u_{k}(\vec{x})$ is
$\mathbf{A}$-holomorphic if and only if there exists the function
$\vec{f'}:\mathbf{A}\rightarrow\mathbf{A}$ such that for all
$k=1,\ldots,n,$ and $\forall \vec{x}\in \mathbf{A}$
\begin{equation}\label{2}
\vec{e}_{k}\vec{f'}(\vec{x})=\lim_{\varepsilon\rightarrow 0}
\frac{\vec{f}(\vec{x}+\varepsilon\vec{e}_{k})-\vec{f}(\vec{x})}{\varepsilon},
\end{equation}
where $\vec{f'}$  doesn't depend on $\vec{e}_{k}.$
\end{theorem}

\begin{proof} Suppose that (\ref{2}) is fulfilled, then it is easily verified that
\begin{equation}\label{3}
\begin{array}{ccc}\vspace*{3mm}
\vec{f'}=\lim_{\varepsilon \rightarrow 0}
\frac{\vec{f}(\vec{x}+\varepsilon\vec{e}_{0})-\vec{f}(\vec{x})}{\varepsilon}=
\sum_{k=0}^{n}\vec{e}_{k}\frac{\partial u_{k}}{\partial
x_{0}},\\\vspace*{3mm}
\vec{e}_{1}\vec{f'}=\lim_{\varepsilon\rightarrow 0}
\frac{\vec{f}(\vec{x}+\varepsilon\vec{e}_{1})-\vec{f}(\vec{x})}{\varepsilon}=
\sum_{k=0}^{n}\vec{e}_{k}\frac{\partial u_{k}}{\partial x_{1}}=
\vec{e}_{1}\sum_{k=0}^{n}\vec{e}_{k}\frac{\partial u_{k}}{\partial
x_{0}},\\\vspace*{3mm} \vdots\\\vspace*{3mm}
\vec{e}_{n}\vec{f'}=\lim_{\varepsilon\rightarrow 0}
\frac{\vec{f}(\vec{x}+\varepsilon\vec{e}_{n})-\vec{f}(\vec{x})}{\varepsilon}=
\sum_{k=0}^{n}\vec{e}_{k}\frac{\partial u_{k}}{\partial x_{n}}=
\vec{e}_{n}\sum_{k=0}^{n}\vec{e}_{k}\frac{\partial u_{k}}{\partial
x_{0}}.
\end{array}
\end{equation}

Consider $\vec{h}=\sum_{k=0}^{n}h_{k}\vec{e}_{k}.$ It follows from
Eqs.(\ref{3}) that
\[
\begin{array}{ccc}\vspace*{3mm}
h_{0}\vec{f'}= h_{0}\sum_{k=0}^{n}\vec{e}_{k}\frac{\partial
u_{k}}{\partial x_{0}},\\\vspace*{3mm} h_{1}\vec{e}_{1}\vec{f'}=
h_{1}\sum_{k=0}^{n}\vec{e}_{k}\frac{\partial u_{k}}{\partial
x_{1}},\\\vspace*{3mm} \vdots\\\vspace*{3mm}
h_{n}\vec{e}_{n}\vec{f'}=
h_{n}\sum_{k=0}^{n}\vec{e}_{k}\frac{\partial u_{k}}{\partial
x_{n}}.
\end{array}
\]

This implies that
\[
\begin{array}{ccc}\vspace*{3mm}
\vec{h}\vec{f'}= h_{0}\sum_{k=0}^{n}\vec{e}_{k}\frac{\partial
u_{k}}{\partial
x_{0}}+h_{1}\sum_{k=0}^{n}\vec{e}_{k}\frac{\partial
u_{k}}{\partial
x_{1}}+\ldots+h_{n}\sum_{k=0}^{n}\vec{e}_{k}\frac{\partial
u_{k}}{\partial x_{n}}=\\\vspace*{3mm}
\lim\limits_{\varepsilon\rightarrow 0}
\frac{\vec{f}(\vec{x}_{0}+\varepsilon\vec{h})-\vec{f}(\vec{x}_{0})}{\varepsilon}.
\end{array}
\]

Furthermore, it follows from Eqs.(\ref{3}) that
\[
\begin{array}{ccc}\vspace*{3mm}
h_{0}\sum_{k=0}^{n}\vec{e}_{k}\frac{\partial u_{k}}{\partial
x_{0}}+h_{1}\sum_{k=0}^{n}\vec{e}_{k}\frac{\partial
u_{k}}{\partial
x_{1}}+\ldots+h_{n}\sum_{k=0}^{n}\vec{e}_{k}\frac{\partial
u_{k}}{\partial x_{n}}=\\\vspace*{3mm}
h_{0}\sum_{k=0}^{n}\vec{e}_{k}\frac{\partial u_{k}}{\partial
x_{0}}+h_{1}\vec{e}_{1}\sum_{k=0}^{n}\vec{e}_{k}\frac{\partial
u_{k}}{\partial
x_{0}}+\ldots+h_{n}\vec{e}_{n}\sum_{k=0}^{n}\vec{e}_{k}\frac{\partial
u_{k}}{\partial x_{0}}.
\end{array}
\]

Therefore, for every $\vec{h}\in A$

$$\vec{h}\sum_{k=0}^{n}\vec{e}_{k}\frac{\partial
u_{k}}{\partial x_{0}}= \lim\limits_{\varepsilon\rightarrow 0}
\frac{\vec{f}(\vec{x}_{0}+\varepsilon\vec{h})-\vec{f}(\vec{x}_{0})}{\varepsilon}
$$
or
\begin{equation}\label{4}
\vec{f'}=\sum_{k=0}^{n}\vec{e}_{k}\frac{\partial u_{k}}{\partial
x_{0}}
\end{equation}
\end{proof}

By using Eq.(\ref{3}), we have

\begin{equation}\label{5}
\begin{array}{ccc}\vspace*{3mm}
\sum_{k=0}^{n}\vec{e}_{k}\frac{\partial u_{k}}{\partial x_{1}}=
\vec{e}_{1}\sum_{k=0}^{n}\vec{e}_{k}\frac{\partial u_{k}}{\partial
x_{0}},\\\vspace*{3mm} \sum_{k=0}^{n}\vec{e}_{k}\frac{\partial
u_{k}}{\partial x_{2}}=
\vec{e}_{2}\sum_{k=0}^{n}\vec{e}_{k}\frac{\partial u_{k}}{\partial
x_{0}},\\\vspace*{3mm} \vdots\\\vspace*{3mm}
\sum_{k=0}^{n}\vec{e}_{k}\frac{\partial u_{k}}{\partial x_{n}}=
\vec{e}_{n}\sum_{k=0}^{n}\vec{e}_{k}\frac{\partial u_{k}}{\partial
x_{0}}.
\end{array}
\end{equation}

Eqs.(\ref{5}) will be called the Cauchy-Riemann type conditions.
It follows from Theorem \ref{A} that if
$\vec{f}(\vec{x})=\sum_{k=0}^{n}\vec{e}_{k}u_{k}(\vec{x})$
satisfies (\ref{5}) then $\vec{f}$ is $\mathbf{A}$-holomorphic.

\begin{theorem}
If $\vec{f}$ is  $\mathbf{A}$-holomorphic and $u_{k}\in
C^{\infty}, k=1,\ldots,n$, then for all $l\geq 1$ there exists
$\vec{f}^{(l)}$, which is $\mathbf{A}$-holomorphic and
$\vec{f}^{(l)}=\sum_{k=0}^{n}\vec{e}_{k}\frac{\partial^{l}
u_{k}}{\partial x_{0}^{l}}.$
\end{theorem}
\begin{proof}It is easy to see that functions $u'_{k}=\frac{\partial
u_{k}}{\partial\vec{x}_{0}}$, $k=1,\ldots,n$, satisfy conditions
(\ref{5}) since $u_{k}\in C^{\infty}$. So $\vec{f}'$ is
$\mathbf{A}$-holomorphic and
$\vec{f}''=\sum_{k=0}^{n}\vec{e}_{k}\frac{\partial^{2}
u_{k}}{\partial x_{0}^{2}}$ (see (\ref{4})). In complete analogy
with this we can show that $\vec{f}^{(l)}$  is
$\mathbf{A}$-holomorphic and
$\vec{f}^{(l)}=\sum_{k=0}^{n}\vec{e}_{k}\frac{\partial^{l}
u_{k}}{\partial x_{0}^{l}}.$
\end{proof}

Let $\vec{f}$ be an $\mathbf{A}$-holomorphic function. For fixed
$\vec{x}, \vec{h}\in \mathbf{A}$ consider the function
$\vec{F}(t)=\vec{f}(\vec{x}+t\vec{h})$. It is easily verified that
$\frac{d^{l}\vec{F}(0)}{dt^{l}}=\vec{f}^{l}(\vec{x})\vec{h}^{l}$.
So the function $\vec{F}(t)$ can be expanded in a Taylor series as
follows
\[
\vec{F}(t)=\sum_{l\geq
0}\frac{1}{l!}\frac{d^{l}\vec{F}(0)}{dt^{l}}t^{l}.
\]
Putting $t=1$, we have
\begin{equation}\label{6}
\vec{f}(\vec{x}+\vec{h})=\vec{f}(\vec{x})+\vec{f'}(\vec{x})\vec{h}+
\frac{1}{2!}\vec{f'}(\vec{x})\vec{h}^{2}+\ldots
\end{equation}

Therefore, every $\mathbf{A}$-holomorphic function (with $u_{k}\in
C^{\infty}$) can be expanded in a Taylor series. In the particular
case where a bicomplex (or hyperbolic) function is
hyperholomorphic, it can be expanded in a Taylor series (\ref{6})
(see \cite{1}, \cite{2}).

\section{Zeros of polynomials in commutative algebras}

Since each  $\mathbf{A}$-holomorphic function can be approximated
by its Taylor polynomial of finite degree, zeros of such functions
might be studied if we can calculate zeros of polynomials. Let
$p_{m}(w)=a_{m}w^{m}+a_{m-1}w^{m-1}+\ldots+a_{0}$ be a polynomial
in the algebra $\mathbf{A}$. Our purpose is to investigate the
structure of the set of zeros of the equation
\begin{equation}\label{7}
 p_{m}(w)=0.
\end{equation}

\begin{theorem}
If $\mathbf{A}$ has $n$ non-trivial idempotents $ i_{1}, i_{2}
,\ldots, i_{n}$ such that $ i_{p}i_{r}=0$ for $p\neq r$, and
$\sum_{l=1}^{n}i_{l}=1$, then Eq. (\ref{7}) can be reduced to the
system of polynomial equations in the field $K$.
\end{theorem}

\begin{proof}
As a preliminary to the proof of the theorem, we shall prove
several auxiliary lemmas.
\begin{lemma}
Idempotents  $ i_{1}, i_{2} ,\ldots, i_{n}$ are linearly
independent vectors.
\end{lemma}
\begin{proof}
Suppose the contrary, then there exist $ k_{1}, k_{2} ,\ldots,
k_{n}\in K$  such that $\sum_{p=1}^{n}|k_{p}|>0$  and
$\sum_{p=1}^{n}k_{p}i_{p}=0$. By using the properties of
idempotents, we have $k_{p}i_{p}=0$  for all $p=1,2,\ldots,n,$ but
this is impossible. Indeed, if $k_{p}i_{p}=0$ for $k_{p}\neq 0,$
then $k_{p}^{-1}(k_{p}i_{p})=i_{p}=0.$
\end{proof}
Denote by $I_{l}=\{ai_{l}|a\in A\}$  the principal ideal generated
by $i_{l}$, $l=1,2,\ldots, n$. It follows from the conditions of
Theorem 3 that the algebra $\mathbf{A}$ can be decomposed in the
direct sum (the Pierce decomposition): $\mathbf{A}=I_{1}\oplus
I_{2}\oplus \ldots \oplus I_{n}$.

\begin{lemma}\label{B}
If  $a\in I_{l}$ then there exists $k\in K$ such that $a=ki_{l}$,
i.e., the ideal $I_{l}$ can be represented in the following form
$I_{l}=\{ki_{l}|k\in K\}.$
\end{lemma}

\begin{proof}
For $a\in I_{l}$ there exists  $b\in \mathbf{A}$  such that
$a=bi_{l}$. Since  $i_{1}, i_{2} ,\ldots, i_{n}$ are linearly
independent, there exist $k_{1}, k_{2} ,\ldots, k_{n}\in K$  such
that $b=\sum_{p=1}^{n}k_{p}i_{p}$. Thus,
$a=bi_{l}=(\sum_{p=1}^{n}k_{p}i_{p})i_{l}=k_{l}i_{l}.$
\end{proof}
Let us consider decompositions
\begin{equation}\label{8}
\begin{array}{ccc}\vspace*{3mm}
a_{r}=a_{r}^{(1)}+\ldots+a_{r}^{(n)}, r=0,1,\ldots
m,\\\vspace*{3mm} w=w_{1}+\ldots +w_{n},
\end{array}
\end{equation}

where $a_{r}^{(p)},w_{p}\in I_{p}.$ Plugging (\ref{8}) into
(\ref{7}), we obtain the following system of polynomial equations
\begin{equation}\label{9}
\begin{array}{ccc}\vspace*{3mm}
a_{m}^{(1)}w_{1}^{m}+a_{m-1}^{(1)}w_{1}^{m-1}+\ldots+a_{0}^{(1)}=0,\\\vspace*{3mm}
a_{m}^{(2)}w_{2}^{m}+a_{m-1}^{(2)}w_{2}^{m-1}+\ldots+a_{0}^{(2)}=0,\\\vspace*{3mm}
\vdots\\\vspace*{3mm}
a_{m}^{(n)}w_{n}^{m}+a_{m-1}^{(n)}w_{n}^{m-1}+\ldots+a_{0}^{(n)}=0.
\end{array}
\end{equation}

It follows from Lemma \ref{B} that $a_{r}^{(s)}=k_{r}^{(s)}i_{s}$,
$w_{s}^{r}=xi_{s}$, where $k_{r}^{(s)},x\in K.$

Therefore, taking $i_{s}$ out of the expression
$a_{m}^{(s)}w_{s}^{m}+a_{m-1}^{(s)}w_{s}^{m-1}+\ldots+a_{0}^{(s)}=0,$
$s=1,\ldots,n,$ the system (\ref{9}) can be reduced to the system
of $n$ polynomial equations in $K$ with coefficients
$k_{r}^{(s)}$.
\end{proof}

\begin{example}
Let $\mathbf{A}$ be the bicomplex algebra, i.e.,
$\mathbf{A}=\{c_{0}+ec_{1}|c_{0},c_{1}\in \mathbb{C}\}$, where
$e^{2}=1$ and $\mathbf{A}$ is commutative. The bicomplex algebra
has two idempotents $i_{1}=\frac{1+e}{2}$ and
$i_{2}=\frac{1-e}{2}$. It is easy to see that $i_{1}i_{2}=0$ and
$i_{1}+i_{2}=1$. Thus, in this case polynomial equation (\ref{7})
can be reduced to the system of two polynomial equations in
$\mathbb{C}$ (see \cite{4}).
\end{example}
\begin{example}
Suppose $\mathbf{A}$ is the commutative algebra of the following
form $\mathbf{A}=\{a_{0}+ea_{1}+fa_{2}+ga_{3}| a_{k}\in
\mathbb{R}\}$, where $e^{2}=f^{2}=g^{2}=1$ and $efg=1$. This
algebra has four idempotents: $i_{1}=\frac{1+e+f+g}{4}$,
$i_{2}=\frac{1-e-f+g}{4}$, $i_{3}=\frac{1+e-f-g}{4}$,
$i_{3}=\frac{1-e+f-g}{4}$. It is easy to see that $i_{k}i_{l}=0$
for $k\neq l$, and $i_{1}+i_{2}+i_{3}+i_{4}=1$. Therefore, in this
case polynomial equation (\ref{7}) can be reduced to the system of
four polynomial equations in $\mathbb{R}$.
\end{example}

\bibliographystyle{amsplain}

\end{document}